\pdfoutput=1
\documentclass[a4paper,11pt]{article}
\usepackage{amsmath,amscd,amsfonts,amssymb,latexsym}

\usepackage{pgf}

\makeatletter
\newtheorem{theorem}{Theorem}
\newtheorem{proposition}[theorem]{Proposition}

\newtheorem{corollary}[theorem]{Corollary}
\newtheorem{remark}[theorem]{Remark}
\newtheorem{definition}[theorem]{Definition}

\newtheorem{example}[theorem]{Example}

\def\proof{\noindent{\bf Proof.\ }}
\def\1{1\kern -.39em1}

\def\qed{\hfill\hbox{$\Box$}}
\def\codim{\mathop{\rm codim}}
\def\Hom{\mathop{\rm Hom}}
\def\Grass{{\mathop{\rm Grass}}}
\def\C{{\mathbb C}}
\def\Z{{\mathbb Z}}

\def\Q{{\mathbb Q}}

\def\O{{\cal O}}
\def\P{{\mathbb P}}
\def\Tt{{\bf T}}
\def\s{\vskip6pt}
\def\bar#1{\overline{#1}}

\begin{document}
\date{ April, 2012}

\title{\bf Equivariant Chern classes\\ and\\ localization theorem}

\author{ Andrzej Weber \\
\small Institute of Mathematics, University of Warsaw\\
\small Banacha 2, 02-097 Warszawa, Poland\\
\small aweber@mimuw.edu.pl }

\maketitle

\begin{abstract} For a complex variety with a torus action we
propose a new method of computing Chern-Schwartz-MacPherson
classes. The method does not apply resolution of singularities. It
is based on the Localization Theorem in equivariant cohomology.\s

This is an extended version of the talk given in Hefei in July
2011.\s


\end{abstract}

Equivariant cohomology is a powerful tool for studying complex
manifolds equipped with a torus action. The Localization Theorem
of Atiyah and Bott and the resulting formula of Berline-Vergne
allow to compute global invariants, for example invariants of
singular subsets, in terms of some data attached to the fixed
points of the action. We will concentrate on the equi\-va\-riant
Chern-Schwartz-MacPherson classes. The global class is determined
by the local contributions coming from the fixed points. On the
other hand, the sum of the local contributions divided by the
Euler classes is equal to zero in an appropriate localization of
equivariant cohomology. Especially for Grassmannians we obtain
interesting formulas with nontrivial relations involving rational
functions. We discuss the issue of positivity: the local
equivariant Chern class may be presented in various ways,
depending on the choice of generating circles of the torus. For
some choices we find that the coefficients of the presentation are
nonnegative. Also the coefficients in the Schur basis are
nonnegative in many examples, but it turns out that not always.

 We
begin with \S\ref{sec1} which contains a review of the results
concerning the equivariant fundamental class of an invariant
subvariety in a smooth $G$-manifold $M$. The first two chapters
are valid for any algebraic group, but further we will consider
only torus
 actions. The equivariant fundamental class lives in the equivariant
cohomology $H^*_G(M)$.  The invariant subvarieties contained in a
vector space $V$, on which the torus $G=T=(\C^*)^n$ acts linearly,
are of particular interest. If the weights of the torus action are
nonnegative then the equivariant fundamental class is a
nonnegative combination of monomials in
$H^*_T(V)=\Q[t_1,t_2,\dots,t_n]$.

In \S\ref{drugi} we discuss the equivariant version of the
Chern-Schwartz-Mac\-Pher\-son class, denoted by $c^T$, which is a
refinement of the equivariant fundamental class. To give the
precise definition, following \cite{Oh}, one has to introduce
equivariant homology. Eventually  we will assume that the variety
is contained in a smooth manifold. By Poincar\'e duality it is
enough for our purpose to consider the equivariant
Chern-Schwartz-MacPherson classes as the elements of equivariant
cohomology of the ambient space.

The main tool of the equivariant cohomology for a torus action is
the Localization Theorem of Atiyah-Bott or Berline-Vergne formula.
It says that the equivariant cohomology class can be read from
certain data concentrated at the fixed points of the action. The
precise formulation of the Localization  Theorem is recalled in
\S\ref{locthe}.

The section \S\ref{rachunki} is a kind of interlude for fun. We
give some examples of calculations based on the Localization
Theorem for torus acting on projective spaces and Grassmannians.
For example we show how the formula for Gysin map for the
Grassmann bundle immediately follows from the Laplace determinant
expansion.

Next, in \S\ref{toryczne} we discuss equivariant
Chern-Schwartz-MacPherson classes of toric varieties. From the
Localization Theorem we deduce that the equi\-va\-riant class of
an orbit $c^T(Tx)$ is equal to the fundamental class of its
closure $[\overline{Tx}]$. Exactly the same formula holds in the
nonequivariant setting, by \cite{BBFK}.

The  section \S\ref{prodwag} is important for whole inductive
procedure of computations of equivariant Chern classes. The key
statement is the following:

\begin{theorem} \label{one} Suppose that $X$ is a $T$-variety,
not necessarily smooth, contained in a $T$-manifold $M$. Let $p\in
X$ be an isolated fixed point. Then the zero degree of the class
$c^T(X)$ restricted to $\{p\}$ is Poincar\'e dual of the product
of weights appearing in the tangent representation
$\Tt_pM$.\end{theorem}

It is also convenient to consider a version of the Localization
Theorem in which we express the global cohomology class by its
restriction to an arbitrary submanifold containing the fixed point
set. The main example is the projective space $\P^n$. The class
which we want to compute is the equivariant
Chern-Schwartz-MacPherson class of the projective cone over a
subvariety in  $\P^{n-1}$. In \S\ref{par-loc}  from partial
localization we deduce the following:

\begin{proposition} Suppose that $X$ is a $T$-invariant cone in a
linear representation $V$ of $T$. Let $h=c_1^T({\O}_{\P(V)}(1))\in
H^2_T(\P(V)) $ be the equivariant Chern class. If
$$c^T({\P(X)})=\Big(\sum_{i=0}^{\dim(V)-1}
b_i(t)h^i\Big)\cap[\P(V)]\in H_*^T( \P(V))$$ for $b_i(t)\in
H^*_T(pt)$ then
 $$c^T(X)=(b_0(t)+e_0)\cap[V]\in H_*^T(V)\,$$
 where $e_0$ is the product of weights appearing in the representation $V$.
\end{proposition}

Even a seemingly trivial application of this result (discussed in
\S\ref{conusy}) is meaningful. If $T=\C^*$ acts by scalar
multiplication on  $\C^n$, then $T$-invariant subvariety is just a
cone in $\C^n$. The characteristic classes of cones were already
considered by Aluffi and Marcolli. In their paper  \cite{AMfe}
there was given a formula for the Chern-Schwartz-MacPherson class
of an open affine cone in $\P^n$. It is not a coincidence, that
their computation agrees with our result about the equivariant
class in $H^*_T(\C^n)$:

\begin{proposition}Suppose $X$ a cone in $\C^n$.
Let $x=c_1(\O_{\P^{n-1}}(1))\in H^2(\P^{n-1})$ and let $t\in
H^2_T(pt)$ be the generator   corresponding to the identity
character. If
$$c_{SM}({\P(X)})=\Big(\sum_{i=0}^{n-1}
a_ix^i\Big)\cap[\P^{n-1}]\in H_*( \P^{n-1})$$
 then
$$c^T({X})= \Big(\sum_{i=0}^{n-1} a_it^i+t^n\Big)\cap[\C^{n}]\in
H_*^T(\C^{n})\,.$$
\end{proposition}
The result follows from the previous one since the equivariant
Chern class $h=c_1^T({\O}_{\P^{n-1}}(1))$ is equal to
$$1\otimes
x-t\otimes 1\in H^*_T(\P^{n-1})=H^*_T(pt)\otimes
H^*(\P^{n-1})\,.$$
 By the product property of equivariant
Chern-Schwartz-MacPherson classes we obtain for free ,,Feynman
rule'' for the polynomial $G_X$ introduced \cite[Lemma
3.10]{AMfe}.

In the next section \S\ref{metoda} we propose a new method of
computing equivariant Chern-Schwartz-MacPherson classes which does
not involve resolution of singularities. It is based on the fact,
that the sum of the equivariant Chern-Schwartz-MacPherson classes
localized at fixed points and divided by Euler classes is equal to
zero, except from the zero degree. A similar observation was
already made by F\'eher and Rim\'anyi in \cite[\S8.1]{FR} for
computation of Thom polynomials. On the other hand the zero degree
of the equivariant Chern class is given by the result of
\S\ref{par-loc} (stated before as Theorem \ref{one}). This way
often we can compute local Chern classes by induction on the depth
of the singularity.

 Our main example in \S\ref{obliczenia} is the determinant variety,
the subset of square matrices $n\times n$ defined by the equation
$\det=0$. We study its compactification, the Schubert variety of
codimension one in $\Grass_n(\C^{2n})$. We discuss computational
problems appearing for that example. The concrete formula for the
equi\-va\-riant Chern class is a huge sum of fractions.
Surprisingly all the difficulties lie in simplifying that
expression. We compute the equiva\-riant Chern class for the
determinant variety for $n\leq 4$. It turns out that it is a
nonnegative combination of monomials with suitable choice of
gene\-ra\-tors of $H^2_T(pt)$. This supports the conjecture of
Aluffi-Mihalcea that the Chern-MacPherson-Schwartz class of the
Schubert varieties are effective. On the other hand for $n=4$ the
local equivariant Chern-Schwartz-MacPherson class expanded in the
Schur basis has negative coefficients in few places. We present
the result of calculations in \S\ref{monstrum}.

The connection of our local formulas with the calculations of
\cite{AlMi} and \cite{Jon} is not clear. The formula for the
global class can be read from the local contributions by Theorem
\ref{wlozenie}. Nevertheless the shape of this relation seems to
be combinatorially nontrivial due to presence of the denominators.

\s I would like to thank Magdalena Zielenkiewicz for correcting
some errors appearing in formulas of the preliminary version of
the paper. Also I would like to thank the Referee for his useful
comments, suggestions and questions that helped a lot to improve
the manuscript.

\section{Equivariant fundamental class}
\label{sec1}

Let $M$ be a complex manifold and $X\subset M$ a closed complex
subvariety. The {\it fundamental class} of $X$, which is the
Poincar\'e dual of the cycle defined by $X$ is denoted by $$[X]\in
H^{2\codim(X)}(M)\,.$$ When the ambient manifold $M$ is
contractible, for example when $M$ is an affine space, there is no
use of $[X]$ since the cohomology of $M$ is trivial. An
interesting situation appears when an algebraic group $G$ acts on
$M$ and $X$ is preserved by the action. In that case there is an
equivariant fundamental class of $X$ which belongs to the
equivariant cohomology of $M$
$$[X]\in H^{2\codim(X)}_G(M).$$
Now even if $M$ is contractible we obtain a remarkable invariant
of the pair $(M,X)$. For contractible $M$ its equivariant
cohomology coincides with the equivariant cohomology of a point
$$H^{2\codim(X)}_G(M)\simeq H^{2\codim(X)}_G(pt)$$
and the cohomology of a point {\it is the ring of characteristic
classes for $G$}. In particular
\begin{itemize}\item if $G=(\C^*)^n$ then $H^*_G(pt)=\Q[t_1,t_2,\dots,t_n]$
\item
if $G=GL_n$ then
$H^*_G(pt)=\Q[\sigma_1,\sigma_2,\dots,\sigma_n]=\Q[t_1,t_2,\dots,t_n]^{\Sigma_n}$
\item in general the ring of characteristic classes coincides with
the invariants of the Weyl group acting on the characteristic
classes for the maximal torus
$$H^*_G(pt)=H^*_T(pt)^W.$$
(We consider here only cohomology with rational
coefficients.)\end{itemize} For a torus $G=T$ we identify
$H^*_T(pt)$ with the polynomial algebra spanned by characters of
$T$, i.e.
$$H^*_T(pt)=\Q[T^\vee]=\bigoplus_{k=0}^\infty Sym^k[T^\vee\otimes\Q]\,.$$
A character $\lambda:T\to\C^*$ corresponds to an element of
$H^2_T(pt)$.

We will briefly recall the construction of equivariant cohomology
in \S\ref{drugi}. The reader can find its basic properties in
\cite{Qu}. For a review of equivariant cohomology in algebraic
geometry see e.g.~\cite{An}.   An extended discussions of
different names for the equivariant fundamental class can be found
in \cite[\S2.1]{BeSe}.

For $G=GL_n$ equivariant cohomology and the equivariant
fundamental classes $[X] \in H^*_{GL_n}(pt)$ has turned out to be
an adequate tool for studying the Thom polynomials of
singularities of maps. Here $X$ is a set of singular jets in the
space of all jets of maps. Its equivariant fundamental class $[X]$
is the universal characteristic class which describes
cohomological properties of singular loci of maps. In the last
decade there appeared a series of papers by Rim\'anyi and his
collaborators (starting from \cite{Ri}) and Kazarian (see
e.g.~\cite{Ka}). Powerful tools allowing effective computations
were developed and some structure theorems were stated. The
geometric approach to equivariant cohomology leads to positivity
results \cite{PW,MPW1,MPW2}. The source of these results is the
following principle:

\begin{theorem} \label{pos0}If
$X\subset \C^N$ is a cone in a polynomial representation of
$GL_n$, then $[X]$ is a nonnegative combination of Schur
functions.\end{theorem}

The examples of polynomial representations are the following: the
natural representation, its tensor products, symmetric products,
exterior products and in general quotients of the sums of tensor
products. The Schur functions constitute a basis of the ring of
characteristic classes
$$H^*_{GL_n}(pt)=H^*(\Grass_n(\C^\infty))$$
corresponding to the decomposition of the infinite Grassmannian
into Schubert cells. For an algebraic treatment of Schur functions
see \cite{McD}.

A version of Theorem \ref{pos0} holds for $G$ being a product of
the general linear groups. We will be interested in torus actions.
Theorem 5 stated in \cite{PW2} reduces to:

\begin{theorem} \label{pos} Let $T=(\C^*)^n$ and let $t_1, t_2,\dots,t_n\in \Hom(T,\C^*)$ be the characters corresponding to the decomposition of $T$ into the
product. Suppose $V=\bigoplus V_\lambda$ is a representation of
$T$ such that each weight $\lambda$ appearing in $V$ is a
nonnegative combination of $t_i$'s. Let $X\subset V$ be a variety
preserved by $T$-action. Then the equivariant fundamental class
$[X]\subset H^*_T(V)=\Q[t_1,t_2,\dots,t_n]$ is a polynomial with
nonnegative coefficients.\end{theorem}

\section{Equivariant Chern class}\label{drugi}
Our goal is to study more delicate invariants of  subvarieties in
representations of algebraic groups, the invariants which are
refinements of the equivariant fundamental class. In most of the
interesting cases the subvarieties to study are singular. Our
first choice is the equivariant version of the
Chern-Schwartz-MacPherson classes. We recall that the usual
Chern-Schwartz-MacPherson classes, introduced in \cite{McP} and
denoted by $c_{SM}$,
 live in homology, they are Poincar\'e duals of the
Chern classes of the tangent bundle when the variety is smooth.
These classes are functorial in a certain sense, and therefore
usually they are computed via resolution of singularities.

The equivariant version of Chern-Schwartz-MacPherson classes was
developed by Ohmoto \cite{Oh}. To define these classes one has to
recall the Borel construction of the equivariant cohomology. Let
$G$ be a topological group. Denote by $EG\to BG=EG/G$ the
universal principal $G$-bundle. This bundle is defined up to
$G$-equivariant homotopy. For a topological $G$-space the
equivariant cohomology is by definition   the cohomology of the
associa\-ted  $X$-bundle
 $EG\times^G X$. Now we apply this construction
to $G$ being an algebraic reductive group and $X$ a complex
algebraic $G$-variety. With the exclusion of the case of the
trivial group, $EG$ does not admit a finite dimensional model.
Instead, $EG$ always has an approximation by algebraic
$G$-varieties, see \cite{To}. For example if $G=\C^*$, then
$EG=\C^\infty-\{0\}$ and $BG=\P^\infty$. It can be approximated by
$U=\C^n-\{0\}$ with $U/G=\P^{n-1}$. In general as the
approximation of $EG$ we take an open set $U$ in a linear
representation $V$ of $G$ satisfying\begin{itemize}
\item $U$ is $G$-invariant,
\item $G$ acts freely on $U$ and the action admits a geometric
quotient,
\item $V-U$ has a sufficiently large codimension in $V$.
\end{itemize}
If we are interested only in the cohomology classes of degrees bounded by
 $d$ we take an approximation with $2\codim(V-U)> d+1$. Then
$$H^k_G(X)=H^k(U\times^G X)$$
for $k\leq 2d$. The equivariant Chern classes of a smooth
$G$-variety coincides with the equivariant Chern class of the
tangent bundle. By Borel construction it is the usual Chern class
of the tangent bundle to fibers of the fibration $EG\times^GX\to
BG$. Using an approximation it can be written as
$$ c^G(X)=p^*c(U/G)^{-1}\cup c(U\times^G X)\in H^*(U\times^G X)\simeq H^*_G(X)\,,$$
where $p:U\times^G X\to U/G$ is the projection. If $X$ is singular
we apply the same formula with obvious modifications. First of all
$U\times^G X$ is singular and we have the homology
Chern-Schwartz-MacPherson class $$c_{SM}(U\times^G X)\in
H^{BM}_*(U\times^G X)\,.$$ (The superscript $BM$ stands for
Borel-Moore homology.) We are forced to use less known equivariant
homology $H^G_*(X)$, \cite{BZ,EdGr2}, which can be defined via
approximation:
$$H^G_k(X)=H^{BM}_{k+2\dim(U/G)}(U\times^G X)$$
for $2n-k<2\codim(V-U)-1$, i.e. for $k>2n-2\codim(V-U)+1$.

\begin{definition}
The equivariant Chern-Schwartz-MacPherson class of $X$ is defined
by the formula
$$ c^G(X)=p^*c(U/G)^{-1}\cap c_{SM}(U\times^G X)\in
H^{BM}_{*+2\dim(U/G)}(U\times^G X)\simeq H^G_*(X)\,.$$ The
definition can be extended to the equivariant constructible
functions on $X$. \end{definition}

Note that $H^G_*(X)$ can have nontrivial negative degrees, but the
equivariant Chern-Schwartz-MacPherson class lives in $H^G_{\geq
0}(X)$.

We will not use the long name {\it equivariant
Chern-Schwartz-MacPherson classes}. Hopefully saying just {\it
equivariant Chern classes} in the context of possibly singular
algebraic $G$-varieties or constructible functions will not lead
to any confusion. Additionally we will always write $c^G(\1_X)$
instead of $c^G(X)$. Later, in \S\S\ref{metoda}-\ref{obliczenia},
where we compute the equivariant Chern classes of a subvariety in
a smooth manifold $M$, for convenience we skip the cap-product
$\cap[M]$ in the notation identifying $H^T_*(M)$ with
$H^{2\dim(M)-*}_T(M)$. \s

The definition of equivariant Chern classes is in fact irrelevant.
All what we need follows from the formal properties.
\begin{itemize}
\item {\bf Normalization:} if $X$ is smooth, then $c^G(\1_X)$ is Poincar\'e dual of the
usual equivariant Chern class of the tangent bundle.
\item {\bf Functoriality:}
for a $G$-constructible function $\alpha$ and a proper $G$-map
$f:X\to Y $ we have $c^G(f_*\alpha)=f_*c^G(\alpha)$.
\item {\bf Product formula:} if $X$ is $G$-variety and $Y$ is
$G'$-variety, then $$c^{G\times G'}(\1_{X\times
Y})=c^{G}(\1_X)\otimes c^{G'}(\1_Y)$$ under the K\"uneth
isomorphism $H^{G\times G'}_*(X\times Y)\simeq H^{G}_*(X)\otimes
H^{G'}_*(Y)$. In particular when $X$ is a trivial $G$-space, then
$$c^G(\1_X)=1\otimes c_{SM}(\1_X)\in H^G_*(X)\simeq H^{-*}_G(pt)\otimes
H_*^{BM}(X)\,.$$

\item {\bf Functoriality with respect to $G$:} Let $\phi:G'\to G$ be a
group homomorphism and $X$ a $G$-space. The induced map
$\phi^*:H^G_*(X)\to H^{G'}_*(X)$ sends $c^G(\1_X)$ to
$c^{G'}(\1_X)$.
\end{itemize}
All five properties easily follow from the corresponding
properties of the usual Chern-Schwartz-MacPherson classes. The
equivariant Chern class carries more information than
nonequivariant Chern-Schwartz-MacPherson class. There is a natural
map $H^G_*(X)\to H_*^{BM}(X)$ which is induced by the inclusion of
the trivial group into $G$. It transports the equivariant Chern
class to the nonequivariant one.

Let us focus on the case when $G=T$ is a torus and $V$ is a
complex linear representation. The equivariant homology $H^T_*(V)$
is a free rank one module over $H_T^*(pt)$ generated by $[V]\in
H^T_{2\dim(V)}(V)$. The action of a character $\lambda\in
H^2_T(pt)$ lowers the degree by 2. By Poincar\'e duality we have
the isomorphisms
$$H^T_{2k}(V)\simeq H^{2(\dim(V)-k)}_T(V)\simeq H^{2(\dim(V)-k)}_T(pt)\simeq Sym^{\dim(V)- k}(T^\vee\otimes\Q)\,.$$
We start with the basic example.

\begin{example} \rm Let $V$ be a complex linear representation of
a torus $T$. Suppose that $V$ decomposes as the sum of the weight
spaces $$V=\bigoplus_\lambda V_\lambda\,.$$ Then
$$c^T(\1_V)=\big(\prod_\lambda(1+\lambda)^{\dim(V_\lambda)}\big)\cap[V]
\in H_*^T(V) \simeq Sym^{\dim(V)-*}(T^\vee\otimes\Q)$$
 and
$$c^T(\1_{\{0\}})=[\{0\}]=\big(\prod_\lambda\lambda^{\dim(V_\lambda)}\big)\cap[V]
\in H_0^T(V) \simeq Sym^{\dim(V)}(T^\vee\otimes\Q)\,.$$ The last
formula follows from covariant functoriality.
\end{example}

Now let us see what the equivariant Chern class means for conical
sets in affine spaces.

\begin{example}\rm\label{AMinv}
Let $T=\C^*$ acts on $\C^n$ by scalar multiplication. Consider a
nonempty
 cone $X\subset \C^n$. We will compute its equivariant Chern
class with respect to the action of $T$. Denote by
$\P(X)\subset\P^{n-1}$ the projectivization of $X$. Let
$h=c_1(\O(1))\in H^2(\P^{n-1})$ and let $t\in H^2_T(pt)$ be the
element corresponding to the identity character. Suppose that
$$c_{SM}(\1_{\P(X)})
=(a_0+a_1 h+\dots +a_{n-1}h^{n-1})\cap [\P^n]\in H_*(\P ^n)\,.$$
We will show in \S\ref{conusy} that the equivariant Chern class of
the cone is equal to
$$c^T(\1_X)=(a_0+a_1t+\dots +a_{n-1}t^{n-1}+t^n)\cap[\C^n]\in H_*^T(\C^n)\,.$$
This formula agrees with computation  of  Aluffi-Marcolli who
calculated the invariant of conical sets defined as the
Chern-Schwartz-MacPherson class of the constructible function
$\1_X$ considered not in $\C^n$ but in $\P^n$.
\end{example}

Now suppose $G=T=(\C^*)^n$ is acting on a vector space $V$ as in
Theorem \ref{pos}. We pose a question: \s

\noindent {\bf Question. }{\it When does $c^T(\1_X)\in
H_*^T(V)\simeq\Q[t_1,t_2,\dots,t_n]$  have nonnegative
coefficients?} \s

This is a special property of $X$ since in general the answer is
negative. If the equivariant Chern classes are effective,
i.e.~represented by an invariant cycle, then the answer is
positive. Also it is easy to find a counterexample: if $T=\C^*$
acts on $V=\C^n$ by scalar multiplication and let $X$ be a cone
over a curve of genus $g>1$ and of degree $d$. Then
$$c^T(\1_X)=([X]+2(1-g)t^{n-1}+t^n)\cap[\C^n]=d[\C^2]+2(1-g)[\C^1]+[\C^0]$$
is a counterexample. On the other hand we have a bunch of positive
examples: local equivariant Chern classes have nonnegative
coefficients for
\begin{itemize}
\item toric singularities (see Corollary \ref{tor}),
\item generic hyperplane arrangements with a small number of hyperplanes
\cite{Alhi},
\item banana Feynman motives \cite{AMba}.
\end{itemize}

\section{Localization theorem}\label{locthe}

For the moment we leave the question of positivity. Our current
goal is to develop a calculus which would allow to compute
equivariant Chern classes avoiding resolution of singularities.
Our main tool is the Localization Theorem for torus action. The
topological setup is the following: suppose the torus $ T=(S^1)^n$
or $(\C^*)^n$ acts on a compact space $M$ (decent enough,
e.g.~equivariant CW-complex). The equivariant cohomology
$H^*_T(M)$ is a module over equivariant cohomology of the point
$$H^*_T(pt)= \Q[t_1,t_2,\dots t_n]\,.$$ The following theorem goes
back to Borel.

\begin{theorem} [\cite{Qu}, \cite{AB}]
The restriction to the fixed set
$$ \iota^*:H^*_T(M)\longrightarrow H^*_T(M^T)$$
becomes an isomorphism after localizing in the multiplicative set
generated by the nontrivial characters
$$S=T^\vee-0\subset H^2_T(pt)\,.$$ \end{theorem}

If $M$ is a manifold, then the inverse of the restriction map is
given by the Atiyah-Bott/Berline-Vergne formula. To explain that
let us fix a notation. We decompose the fixed point set into
components $M^T= \bigsqcup_{\alpha\in A} M_\alpha$. Each
$M_\alpha$ is a manifold and denote by $e_\alpha\in
H^*_T(M_\alpha)=H^*(M_\alpha)\otimes \Q[t_1,t_2,\dots t_n]$ the
equivariant Euler class of the normal bundle. The following map is
the inverse of the restriction to the fixed points
\begin{equation} \label{wlozenie1}\begin{matrix}
S^{-1}H^*_T(M^T)&=&\bigoplus_{\alpha\in A}
S^{-1}H^*_T(M_\alpha)&\stackrel{\simeq}{\longrightarrow}&
S^{-1}H^*_T(M) \cr \phantom{.}\cr &&\{ x_\alpha\}_{\alpha\in A}&
\mapsto &\sum_{\alpha\in A} {\iota_\alpha}_*\left(\frac{x_\alpha}{
e_\alpha}\right)\,,\end{matrix}\end{equation} where $
\iota_\alpha:M_\alpha\hookrightarrow M$ is the inclusion. The key
point in the  formula (\ref{wlozenie1}) is that the Euler class $
e_ \alpha$ is {\it invertible} in $ S^{-1}H^*_T(M_ \alpha)$.

\begin{remark}\rm Note that if $M$ is a smooth compact algebraic variety and
the action of the torus is algebraic, then $H^*_T(M)$ is a free
module over $H^*_T(pt)$, so we do not kill any class inverting
nontrivial characters.\end{remark}

I other words we can state the theorem:

\begin{theorem}[\cite{AB,EdGr}]\label{wlozenie} Let $M$ be an algebraic variety with algebraic torus action, then with the previous notation
$$x=\sum_{\alpha\in A} {\iota_\alpha}_*\left(\frac{x_{|M_\alpha}}{
e_\alpha}\right)\in H^*_T(M)\,.$$
\end{theorem}

Therefore we can say that $x$ is a sum of local contributions.
Although one has to understand that this statement is a bit
misleading. In fact it is not possible to extract individual
summands in $H^*_T(M)$. This can be done only in the localized
ring. How a single fixed point component contributes to the global
class is obscured by the weights of the tangent representation.

Furthermore consider the push-forward, i.e.~the integration along
$ M$
$$ p_*=\int_M:H^*_T(M)\to H^{*-2\dim(M)}_T(pt)$$
where $ p:M\to pt$ is the constant map. Another form of the
Localization Theorem allows to express the integration along $M$
by integrations along components of the fixed point set.

\begin{theorem} [Berline-Vergne \cite{BV}] For $x\in H^*_T(M)$ the integral can be computed by
summation of local contributions
\begin{equation}\int_M x=\sum_ {\alpha\in A} \int_{M_ \alpha}
\frac{x_{|M_ \alpha}}{ e_
\alpha}\,.\label{integral}\end{equation}\end{theorem}

In particular, when the fixed point set is discrete
$M^T=\{p_0,p_1,\dots p_n\}$ then the Euler class is the product of
weights
$$e_{p}=\prod_ {\lambda\in\Lambda} \lambda^{\dim(V_\lambda)}
\in H_T^*(\{p\})\in\Q[t_1,t_2,\dots,t_n]\,,$$ provided that
$\Tt_{p}M$, the tangent space at $p$ is the sum of weight spaces
$$\Tt_{p}M= \bigoplus_{\lambda\in\Lambda} V_\lambda\,.$$
The integral along $M$ is equal to the sum of fractions:
$$ \int_M a=\sum_{p\in M^T}
\frac{a_{|p}}{ e_{p}}\,.$$

\begin{remark}
\rm The Berline-Vergne formula (\ref{integral}) can be formulated
for singular spaces embedded into a smooth manifold. The local
factor $\frac{1}{e_\alpha}$ is replaced by
$\frac{[X]_{|M_\alpha}}{e_\alpha}$, see \cite{EdGr,BeSe}. There is
a generalization of the Theorem (\ref{wlozenie}) for equivariant
homology (or Chow groups) of singular spaces, but one needs an
additional assumption allowing to define $\iota^*$,
\cite[Proposition 6]{EdGr}.
\end{remark}

\section{Some calculi of rational functions}\label{rachunki}

Before examining equivariant Chern classes of Schubert varieties
let us look closer at some computations based on the Localization
Theorem for Grassmannians. Let us start with the projective space
$M=\P^n$ with the standard torus $T=(\C^*)^{n+1}$ action. The
fixed point set is discrete and consists of coordinate lines
$$M^T=\{p_0,p_1,\dots p_n\}\,.$$ The
tangent space at the point $p_k$ decomposes into one dimensional
representations:
$$\Tt_{p_k}M= \bigoplus_{\ell\not= k}\C_{t_\ell-t_k}\,.$$
The Euler class is equal to
$$e_{p_k}=\prod_{\ell\not= k} (t_\ell-t_k)\,.$$
Let us integrate powers of $c_1:=c_1(\O(1))$. Of course
$$\int_{ \P^n} c_1^m=\left\{\begin{matrix}0\qquad {\rm for}\quad m<n\cr 1\qquad {\rm for}\quad m=n\end{matrix}\right.$$
Applying Berline-Vergne (\ref{integral}) formula we get the
identity
\begin{equation} \sum_{k=0}^n \frac{(-t_k)^m}{ \prod_{\ell\not=k}
(t_\ell-t_k)}=\left\{\begin{matrix}0& {\rm for}\quad m<n\hfill\cr
1& {\rm for}\quad
m=n\,,\hfill\end{matrix}\right.\label{rzutowa}\end{equation} which
is not obvious at the first sight. For example we encourage the
reader to compute by hand the sum \s\s\noindent
$\frac{t_0^2}{(t_1-t_0)(t_2-t_0)(t_3-t_0)}+\frac{t_1^2}{(t_0-t_1)(t_2-t_1)(t_3-t_1)}+
\frac{t_2^2}{(t_0-t_2)(t_1-t_2)(t_3-t_2)}+\frac{t_3^2}{(t_0-t_3)(t_1-t_3)(t_2-t_3)}.$
\s\s \noindent This is exactly the expression (\ref{rzutowa}) for
$m=2$, $n=3$. Replacing $$t_0=0\,,\quad t_1=1\,,\quad t_2=2\,,\;
\dots\, ,\;t_n=n$$ (i.e. specializing to a subtorus) the sum
(\ref{rzutowa}) is equal to
$$\sum_{k=0}^n \frac{(-1)^{m+k}k^m}{ k!(n-k)!}$$
Multiplying by $n!$ we obtain
$$ \sum_{k=0}^n \binom{n}{k}(-1)^{m+k}k^m=\left\{\begin{matrix}0&{\rm for}\quad m<n\hfill\cr n!& {\rm for}\quad m=n\,.\hfill\end{matrix}\right.$$
which is a good exercise for students.

The integral of higher powers of $c_1$ is even more interesting:
Let us see what we do get for $m>n$? For example $n=2$, $m=4$ we
have
$$\frac{t_0^4}{(t_1-t_0)(t_2-t_0)}+\frac{t_1^4}{(t_0-t_1)(t_2-t_1)}+\frac{t_2^4}{(t_0-t_2)(t_1-t_2)}$$
It takes some time to check that the sum is equal to
$$t_0^2 +t_1^2+ t_2^2 + t_0 t_1 + t_0 t_2 + t_1 t_2,.$$
In terms of the elementary symmetric functions it is equal to
$$\sigma_1^2-\sigma_2\,.$$
\begin{proposition}\label{pushf} In general
$$(-1)^k\int_{ \P^n} c_1^{n+k}$$
is equal to the Schur function $S_k$ (which corresponds to the
Segre class of vector bundles).\end{proposition}

\proof By Jacobi-Trudy formula (which is the definition of the
Schur function)
$$S_k(t_0,t_1,\dots,t_n)=\frac{
\left| \begin{matrix} t_0^{n+k} & t_0^{n-1} & t_0^{n-2} &\dots &
t_0^1&1 \cr t_1^{n+k} & t_1^{n-1} & t_1^{n-2} &\dots & t_1^1&1 \cr
t_2^{n+k} & t_2^{n-1} & t_2^{n-2} &\dots & t_2^1&1 \cr
\vdots&\vdots&\vdots&&\vdots&\vdots \cr t_n^{n+k} & t_n^{n-1} &
t_0^{n-2} &\dots & t_n ^1&1 \cr
\end{matrix}\right|}{\prod_{i<j}(t_i-t_j)}\,.$$
To prove the proposition it is enough to use Laplace expansion
with respect to the first column and watch carefully the
signs.\qed\s

\begin{remark} \rm It is wiser to use the dual Grassmannian of hyperplanes, then one gets rid of the factor $(-1)^k$. The general formula with positive signs for Grassmannians is given by Theorem \ref{jlp}.\end{remark}

We will have a look now at the calculus on $\Grass_m(\C^n)$. The
fixed point set consists of coordinate subspaces:
$$\Grass_m(\C^n)^T=\{p_\lambda\,:\,\lambda=(\lambda_1<\lambda_2<\dots<\lambda_m), \;1\leq\lambda_1,\;\lambda_m\leq n\}$$
The tangent space at the fixed point $p_\lambda$ decomposes into
distinct line representations of $T$:
$$\Tt_{p_\lambda}\Grass_m(\C^n)=\bigoplus_{k\in \lambda, \ell\not\in
\lambda}\C_{t_\ell-t_k}\,.$$ The Euler class is equal to
$$ e_{p_\lambda}=\prod_{k\in \lambda,\;\ell\not\in
\lambda} (t_\ell-t_k)\,.$$ Let us integrate a characteristic class
of the tautological bundle ${\cal R}_m$. Suppose that the class
$\phi({\cal R}_m)$ is given by a symmetric polynomial in Chern
roots $W(x_1,x_2,\dots,x_m)$. Then
$$\int_{\Grass_m(\C^n)}\phi({\cal R}_m)
=\sum_\lambda\frac{W(t_i:i\in \lambda)}{\prod_{k\in
\lambda,\;\ell\not\in \lambda} (t_\ell-t_k)}$$ It looks like a
rational function, but we obtain a polynomial in $t_i$'s of degree
$deg(W)-\dim(\Grass_m(\C^n))$.  This expression can be written as
the iterated residue\footnote{More general formulas were found
recently by Magdalena Zielenkiewicz for Grassmannians of all
classical groups.}
\begin{equation}\frac1{m!}Res_{z_1=\infty}Res_{z_2=\infty}\dots
Res_{z_m=\infty}\frac {W(z_1,z_2,\dots
z_m)\prod_{i\not=j}(z_i-z_j)}{\prod_{i=1}^n\prod_{j=1}^m(t_i-z_j)}\,,\label{residua}\end{equation}
see \cite{Be}. Of course if $deg(W)<\dim(\Grass_m(\C^n))=(n-m)m$,
then
$$ \sum_\lambda \frac{W(t_i:i\in\lambda)}{\prod_{k\in \lambda,\;\ell\not\in \lambda}
(t_\ell-t_k)}=0$$ If $deg(W)=\dim(\Grass_m(\C^n))$, then we get a
constant. For example for
$W=c_1^{\dim(\Grass_m(\C^n))}=\big(-(x_1+x_2+\dots+x_m)\big)^{(n-m)m}$
we obtain the degree of the Pl\"ucker embedding
$\Grass_m(\C^n)\subset\P(S^m(\C^n))$ (or the volume of
$\Grass_m(\C^n)$). According to Hook Formula \cite[\S4.3]{FY}
$$deg(\Grass_m(\C^n))=\frac{(m(n-m))!}{\prod_{(i,j)\in\lambda} h(i,j)}\,,$$
where $h(i,j)$ denotes the length of the hook with vertex at
$(i,j)\in\lambda$ contained in the rectangle $m\times (n-m)$. For
$\Grass_3(\C^7)$ the hook lengths are the following
$$\begin{array}{|c|c|c|c|c|}\hline 6 & 5 & 4 & 3\cr\hline
5 & 4 & 3 & 2\cr\hline 4 & 3 & 2 & 1\cr\hline\end{array}$$ Hence
the degree is equal to
$$\frac {12!}{6 \cdot 5 \cdot 4 \cdot 3 \cdot
5 \cdot 4 \cdot 3 \cdot 2 \cdot 4 \cdot 3 \cdot 2 \cdot 1 }= 462.
$$
It would be interesting to find an immediate connection of the
Hook formula and the residue method given by the formula
(\ref{residua}). \s

Let us now formulate a generalization of Proposition \ref{pushf}.
For a partition $I=(i_1\geq i_2\geq i_n)$ the Schur function is
defined by Jacobi-Trudy formula \cite[\S I.3]{McD}
$$S_I(t_1,t_2,\dots, t_n)=\frac{
\left| \begin{matrix} t_1^{n-1+i_1} & t_1^{n-2+i_2} &
t_1^{n-3+i_3} &\dots & t_1^{1+i_{n-1}}&t_1^{i_n} \cr t_2^{n-1+i_1}
& t_2^{n-2+i_2} & t_2^{n-3+i_3} &\dots & t_2^{1+i_{n-1}}&t_2^{i_n}
\cr t_3^{n-1+i_1} & t_3^{n-2+i_2} & t_3^{n-3+i_3} &\dots &
t_3^{1+i_{n-1}}&t_3^{i_n} \cr \vdots&\vdots&\vdots&&\vdots&\vdots
\cr t_n^{n-1+i_1} & t_n^{n-2+i_2} & t_n^{n-3+i_3} &\dots &
t_n^{1+i_{n-1}}&t_n^{i_n} \cr
\end{matrix}\right|}{\prod_{1\leq i<j\leq n}(t_i-t_j)}\,.$$
The definition of Schur function is extended to characteristic
classes of vector bundles. Expanding the determinant with respect
to the first block column containing $m\times m$ minors we find
the formula for push-forward:
\begin{theorem} \label{jlp} Consider the quotient bundle $\cal Q$ and the
tautological bundle $\cal R$ over $\Grass_m(\C^{n})$. Let $J=(j_1\geq j_2\geq\dots\geq j_{n-m})$ and $K=(k_1\geq k_2\geq\dots\geq k_m)$
be partitions. Suppose $j_{n-m}-m\geq k_1$. Then
$$\int_{\Grass_m(\C^{n})}S_J({\cal Q})S_K({\cal R})=S_I(t_1,t_2,\dots,t_n)\,,$$
where $I=(j_1-m\geq j_2-m\geq\dots\geq j_{n-m}-m\geq k_1\geq
k_2\geq\dots\geq k_m)$.
\end{theorem}
A suitable modifications of Theorem \ref{jlp} can be easily
formulated for the partitions not satisfying the inequality
$j_{n-m}-m\geq k_1$. The integral is equal up a sign to the Schur
function for another partition or it is zero. By the splitting
principle Theorem \ref{jlp} implies the corresponding statement
for Grassmannian bundles over any base, not necessarily over the
classifying space $BT$. This way we obtain a proof of the Gysin
homomorphism formula \cite{JLP}, \cite[\S4.1]{FP}. \s

The equivariant Schubert calculus was studied by a number of
authors: Knutson--Tao \cite{KnTa}, Laksov--Thorup \cite{La-Th},
Gatto--Santiago \cite{Ga-Sa} and others. Some formulas can be
obtained by taking residue at infinity \cite{Be,BeSe}. Concluding
this section I would like to say that it seems that still the
calculus of {\it rational} symmetric functions is not developed
enough. In \S\ref{detvar} we will present a method of computation
of equivariant Chern classes of Schubert varieties. Unfortunately
I do not know (maybe except Theorem \ref{jlp}) a tool which would
allow us to simplify the expressions which appear in computation.

\section{Toric varieties}\label{toryczne}
We keep in mind that our purpose is to compute equivariant Chern
classes. From the Localization Theorem it follows that equivariant
Chern classes are determined by local equivariant Chern classes
belonging to the homologies of the components of $M^T$. In the
beginning let us consider the toric varieties, which are quite
easy, but unfortunately not very general from our point of view.
\s

\begin{theorem} Let $X$ be a toric variety.
Consider the cycle $\Xi_X$ which is equal to the sum of the
closures of orbits. Then $\Xi_X$ represents the equivariant Chern
class $c^T(\1_X)\in H^T_*(X)$.\end{theorem}

\proof First we consider the case when $X$ is a smooth toric
variety. If $X=\C^1$ with the standard action of $T=\C^*$ then,
indeed, the equivariant Chern class is equal $[\C]+[0]=[\Xi_X]$.
 By Whitney formula and the product property of sets the statement
 holds for $X=\C^n$ with the standard action
 of $T=(\C^*)^n$.
 Every smooth toric variety locally looks like $\C^n$ with the standard action of the torus, therefore the equation $c^T(1_X)=[\Xi_X]$ holds locally, i.e.~after restriction to
each fixed point.
 Let $X$ be a complete smooth toric variety. Then $H_*^T(X)$
is free over $H^*_T(pt)$. By the Localization Theorem
$c^T(\1_X)=\Xi_X$ holds globally.
 The noncomplete case follows since any smooth toric variety can be compactified
equivariantly.

 The singular case can be deduced as usual by
functoriality. One sees that for smooth toric varieties the
equivariant Chern class of the constructible function supported by
a single orbit is exactly the fundamental class of the closure of
that orbit without boundary cycles. The equality is preserved by
the push-forward. \qed \s

Note that the theorem holds in equivariant homology of $X$ and we
do not have to use any embedding into a smooth manifold. The
non-equivariant case was proven by Ehlers and
Barthel-Brasselet-Fieseler \cite{BBFK} and it also follows
immediately from \cite{AlCh}.

The cycle representing the equivariant Chern class of a toric
variety is effective. Therefore for the embedded case by Theorem
\ref{pos} we have the corollary:

\begin{corollary}\label{tor} Let $V$ be a representation of $T$.
Suppose an affine $T$-variety $X$ (possibly singular) is embedded
equivariantly into $V$. If the weights of the torus acting on $V$
are nonnegative then the coefficients of $c^T(\1_X)\in
H^*_T(V)=\Q[t_1,t_2,\dots,t_n]$ are nonnegative.\end{corollary}
The situation described in the Corollary \ref{tor} appears when
$$X=X_\sigma=Spec(\C(\sigma^\vee\cap N))$$ is presented in the
usual way: the embedding into $$V=Spec(\C[x_1,x_2,\dots,x_n])$$ is
given by a choice of the generators of the semigroup
$\sigma^\vee\cap N$, see \cite[\S1.3]{Futor}. \s

\begin{remark}\rm All the singularities of the Schubert varieties in
Grassmannians of planes $Grass_2(\C^n)$ are toric. Therefore the
local equivariant Chern classes are nonnegative combinations of
monomials for a suitable choice of a basis of $H^2_T(pt)$.
\end{remark}

The global positivity of Chern classes of Schubert varieties in
$Grass_2(\C^n)$ seems not to follow automatically. Except from the
case of projective spaces only the Schubert varieties with
isolated singularities (the partitions $(n-3,k)$ for $k\leq n-4$,
according to the standard convention) are toric. Nevertheless it
was shown in \cite[\S4.3]{AlMi} that the  nonequivariant  Chern
classes of Schubert varieties in $Grass_2(\C^n)$ are indeed
effective.

\section{Equivariant Chern class of degree zero}\label{prodwag}
The following Theorem \ref{eu} is the key to the inductive
procedure for computing equivariant Chern classes. The theorem
says that the degree zero component of the equivariant Chern class
localized at a fixed point does not depend seriously on the set
itself, but only on wether the point belongs to the set or not.

\begin{theorem} \label{eu} Suppose that $X$ is a $T$-variety,
not necessarily smooth, contained in a $T$-manifold $M$. Let $p\in
X$ be an isolated fixed point. Then the  degree zero component of
the class $c^T(\1_X)$ restricted to $\{p\}$ is Poincar\'e dual to
the product of weights appearing in the tangent representation
$\Tt_pM$
$$\left(c^T(\1_X)_{(0)}\right)_{|p}=e_p\cap[p]\,.$$
\end{theorem}
 By additivity of equivariant Chern classes
it follows that if $p\not\in X$ then $(c^T(\1_X)_{(0)})_{|p})=0$.

The core of the proof is the basic equation of Euler
characteristics
$$\chi(X)=\chi(X^T)\,.$$
Nevertheless the argument demands some formal manipulations.
First of all we note the following fact.
 \begin{proposition}\label{euclas} Let $N$ be a complete
 manifold with a torus action. Let us
decompose the fixed point set $N^T=\sqcup_{\alpha\in A}N_\alpha$
into connected components. Let $i_\alpha:N_\alpha\to N$ be the
inclusion.  The  equivariant cohomology top Chern class of  $N$ is
equal to the sum
\begin{equation}c^T_{top}(N)=\sum_{\alpha\in
A}(i_\alpha)_*\left(c_{top}( N _\alpha)\right)\in H^{2\dim(N)}_T(
N )\,.\end{equation} Dually we have
\begin{equation}c^T(\1_N)_{(0)}=\sum_{\alpha\in
A}(i_\alpha)_*\left(c_{SM}( \1_{N_\alpha})_{(0)}\right)\in H_{0}^T(
N )\,.\end{equation}
\end{proposition} \proof The proof is the
straightforward application of the Theorem (\ref{wlozenie}) since
$$i_\alpha^*(c_{top}^T(N))=e_\alpha\cdot c_{top}(N_\alpha)\in
H^*_T(N_\alpha)= H^*_T(pt)\otimes H^*(N_\alpha)\,.$$\qed\s

\noindent {\bf Proof of Theorem \ref{eu}.} Denote by
$\iota_p:\{p\}\to M$ the inclusion of the point. We will argue
that for any equivariant constructible function $\alpha:M\to\Z$
the equality holds
$$\iota_p^*(c^T(\alpha))=\alpha(p)\,e_p\cap[p]\in H_{-2\dim(M)}^T(\{p\})\,.$$
It is enough to show that statement for $M$ complete and the
constructible function of the shape $\alpha=f_*(1\!\!1_N)$ for an
equivariant map $f:N\to M$ from a smooth complete variety $N$. (We
can assume that $N$ is smooth by the usual argument which is
available thanks to equivariant completion \cite{Sum} and
equivariant resolution of singularities \cite{BM}.) It remains to
prove that
\begin{equation}\iota_p^*f_*(c^T(\1_N)_{(0)})=
\chi(f^{-1}(p))\,e_p\cap[p]\in
H_{-2\dim(M)}^T(\{p\})\,.\label{obrazeu}\end{equation} Let
$i_\alpha$ be as in Proposition \ref{euclas} and
$f_\alpha=f\,i_\alpha:N_\alpha\to M$. We compute the push-forward
of the zero degree component:
 \begin{align}f_*c^T(\1_N
)_{(0)} &=\sum_{\alpha\in A}f_*(i_\alpha)_*\big(c_{SM}( \1_{N
_\alpha})_{(0)}\big)\cr &=\sum_{\alpha\in
A}(f_\alpha)_*\left(c_{SM}(\1_{ N _\alpha})_{(0)}\right) \in
H^T_0(M)\,.
\end{align}
Let $B\subset A$ be the set of components of $N^T$ which are
mapped to $p$. Then
\begin{align}\iota_p^*\big(f_*c^T( \1_N )_{(0)}\big)
&=\iota_p^*\big(\sum_{\beta\in B}(f_\beta)_*c(\1_{ N
_\beta})_{(0)}\big)\cr &=\iota_p^*\sum_{\beta\in B}\chi( N
_\beta)\,[p]\cr &=\sum_{\beta\in B}\chi( N _\beta)\,e_p\cap[p]\,.
\end{align}
We conclude that the equation (\ref{obrazeu}) holds because
$\chi(f^{-1}(p))=\chi(f^{-1}(p)^T)$ and
$f^{-1}(p)^T=\bigsqcup_{\beta\in B}N_\beta$. \qed

\section{Partial localization}\label{par-loc}

There exists the following modification of the localization
formula: we can replace $M^T$ by any invariant submanifold or even
arbitrary invariant subset) $Y$ containing the fixed point set
$M^T$. Then the restriction map
$$H^*_T(M)\to H^*_T(Y)$$
becomes an isomorphism after inversion of nontrivial characters
$S$. Also the Berline-Vergne formula holds, but it makes sense
only for $Y$ being a submanifold. Suppose that $Y=Y_1\sqcup
\{p\}$. It follows that for any $x\in H^*_T(M)$ we have
\begin{equation}\frac{x_{|p}}{e_p}+\int_{Y_1}\frac{x_{|Y_1}}{e_{Y_1}}=0\label{locmod}\end{equation}
for degree smaller than $\dim(M)$. We will apply this formula for
Poincar\'e dual of $c^T(\1_X)$. The integral of the zero degree
Chern-Schwartz-MacPherson class (which corresponds to the top
degree of the cohomology class) is equal to the Euler
characteristic and the same holds for the equivariant Chern class
by the commutativity of the diagram:
$$\begin{matrix}c^T(\1_X)_{(0)}&\in& H^T_0(M)&\longrightarrow&H_0(M)&\ni&c(X)_{(0)}\cr
\phantom.\cr &&\downarrow&&\downarrow\cr \phantom.\cr \int_M
c^T(\1_X)_{(0)}&\in&
H^T_0(pt)&\stackrel{\simeq}{\longrightarrow}&H_{0}(pt)&\ni&\chi(X)\,.\cr
\end{matrix}$$
We apply the partial localization and we find that
\begin{equation}\frac{(c^T(X)_{(0)})_{|p}}{e_p}+
\int_{Y_1}\frac{(c^T(X)_{(0)})_{|Y_1}}{e_{Y_1}}=\chi(X)\,.\label{locmodtop}\end{equation}
Here $e_{Y_1}$ is the equivariant Euler class of the normal bundle
of $Y_1$. (Of course it may be of different degrees over distinct
components of $Y_1$.)

\begin{example} \rm The partial localization allows us to compute the equivariant Chern class of the affine
cone over a projective variety. Suppose $T$ acts on $\C^n$ with
nonzero weights
$$w_1,w_2,\dots,w_n\,.$$
First recall that the equivariant cohomology ring of $\P^{n-1}$ is
the quotient of the polynomial algebra
$$H^*_T(pt)[h]=\Z[t_1,t_2,\dots,t_n,h]$$
by the relation
$$\prod_{i=1}^n(h+w_i)=0. $$
Using the elementary symmetric functions $\sigma_i$ the relation
takes form
\begin{equation}\sum_{i=0}^n\sigma_i(w_\bullet)h^{n-i}=0\,.
 \label{relacja}\end{equation}
Let $X\subset \C^n$ be a nonempty $T$-invariant cone and
$\P(X)\subset \P^{n-1}$ its projectivization. We consider $X=\bar
X-\P(X)$ as a constructible set in $\P^n$ and we will compute its
equivariant Chern class in $H^*_T(\P^n)$.  In this example we skip
the Poincar\'e duals in the notation. Denote by
$\iota:\P^{n-1}\to\P^n$ the inclusion. The equivariant Chern class
of $X$ restricted to $\P^{n-1}$ is equal to
 $$\begin{matrix}\iota^*c^T(\1_X)&=\iota^*c^T(\1_{\bar
X})-\iota^*\iota_*c^T(\1_{\P(X)})\hfill\\ \\&= (1+h)\cdot
c^T(\1_{\P(X)})-h\cdot c^T(\1_{\P(X)})\hfill\\
\\&=c^T(\1_{\P(X)})\,.\hfill\end{matrix}$$
Suppose that the equivariant Chern class of $\P(X)$ is written as
$$c^T(\1_{\P(X)})=\sum_{i=0}^{n-1} b_i(t)h^i\in H^*_T( \P^{n-1})$$
for some polynomials $b_i(t)\in H^*_T(pt)$ of degree $\leq n-1-i$.
To compute the local equivariant Chern class at $0$ we will apply
the formulas (\ref{locmod}) and (\ref{locmodtop}) to $M=\P^n$,
$Y=\{0\}\cup\P^{n-1}$ and $Y_1=\P^{n-1}$. We compute
$$\int_{\P^{n-1}}\frac{c^T(\1_X)}{e_{Y_1}}=\int_{\P^{n-1}}\sum_{i=0}^{n-1} b_i(t)h^{i-1}\,.$$
Except from $i=0$ the summands are integral (belong to $H^*_T(
\P^{n-1})$) and they are of  degree smaller than $n-1$. Therefore
$$\int_{\P^{n-1}}\frac{c^T(\1_X)}{e_{Y_1}}=\int_{\P^{n-1}}\frac{b_0(t)}h\,.$$
An easy calculation using (\ref{relacja}) shows that the inverse
Euler class of the normal bundle to $\P^{n-1}$ is equal to
$$h^{-1}=-\sum_{i=1}^n
\frac{\sigma_{n-i}(w_ \bullet)}{\sigma_{n}(w_
\bullet)}h^{i-1}\,.$$ Hence
$$\int_{\P^{n-1}}\frac{c^T(\1_X)}{e_{Y_1}}=-\int_{\P^{n-1}}b_0(t)\sum_{i=1}^n
\frac{\sigma_{n-i}(w_ \bullet)}{\sigma_{n}(w_
\bullet)}h^{i-1}=-\frac{b_0(t)}{\sigma_{n}(w_ \bullet)}$$ By the
formulas (\ref{locmod}) and (\ref{locmodtop}) and since
$\sigma_{n}(w_ \bullet)=e_p$ we find that
$$\frac{c^T(\1_X)_{|p}}{e_p}-\frac{b_0(t)}{e_p}=\chi(X)=1\,.$$
Therefore
$$c^T(\1_X)_{|p}=b_0(t)+e_0\in H^*_T(\{p\})\,$$
\end{example}

We obtain the following result:

\begin{proposition}\label{stozek}Suppose that $X$ is a nonempty $T$-invariant cone in a
linear representation $V$ of $T$. Let $h=c_1^T({\O}_{\P(V)}(1))$
be the equivariant Chern class. If
\begin{align}c^T(\1_{\P(X)})=\Big(\sum_{i=0}^{\dim(V)-1}
b_i(t)h^i\Big)\cap[\P(V)]\in H_*^T(
\P(V))\label{zdzial}\end{align} then
 $$c^T(\1_X)=(b_0(t)+e_0)\cap[V]\in H_*^T(V)\,$$
 where $e_0$ is the Euler class of the representation $V$.
\end{proposition}
\proof First note, that restriction $H^*_T(V)\to  H^*_T(pt)$ is an
isomorphism. We apply the calculation of the previous example. The
degree of $b_0$ is at most $\dim(V)-1$, therefore it does not
interfere with $e_p$, which is homogeneous of degree
$\dim(V)$.\qed

\section{Conical sets in an affine space}
\label{conusy}

 We come back to the  Example \ref{AMinv} of \S\ref{drugi} which was the starting
point of our interest in equivariant Chern classes. In \cite{AMfe}
there was defined an invariant of a conical set $X\subset \C^n$.
It is equal to the Chern-Schwartz-MacPherson class of $X$
considered as a constructible set in $\P^n$. This Chern class
$$c_{SM}(\1_X)\in H^*(\P^n)$$ is expressed via the Chern class of
the projectivization. The calculation is based on the following
formula:

\begin{proposition}[{\cite[Prop 5.2]{AMba}}] \label{AlMafro}
Let  $X\subset\C^n$ be a nonempty conical set. Let $\bar X=X\cup
\P(X)$ be the closure of $X$ in $\P^n$. Let
$x=c_1(\O_{\P^{n-1}}(1))$ and $\tilde x=c_1(\O_{\P^{n}}(1))$.
Suppose that
\begin{align}\label{zwykla0}c_{SM}(\1_{\P(X)})=\Big(\sum_{i=0}^{n-1}
a_ix^i\Big)\cap[\P^{n-1}]\in H_*( \P^{n-1})\end{align}
 then
\begin{align}\label{zwykla}c_{SM}(\1_{\bar X})=\Big((1+\tilde x)\big(\sum_{i=0}^{n-1}
a_i\tilde x^i\big)+\tilde x^n\Big)\cap[\P^{n}]\in H_*(
\P^{n})\end{align}
\end{proposition}
It follows that
\begin{align}c_{SM}(\1_X)=\Big(\sum_{i=0}^{n-1}
a_i\tilde x^i+\tilde x^n\Big)\cap[\P^{n}]\in
H_*(\P^{n})\end{align}

It seems natural to look at the conical sets from the point of
view of equivariant cohomology. Let $T=\C^*$ acts on $\C^n$ by
scalar multiplication.
\begin{proposition}\label{antal}Under assumption of Proposition \ref{AlMafro}
\begin{align}c^T(\1_{ X})=
\Big(\sum_{i=0}^{n-1} a_it^i+t^n\Big)\cap[\C^{n}]\in
H_*^T(\C^{n})\end{align}
\end{proposition}

\proof Assume that the usual, nonequivariant Chern class of
$\P(X)$ satisfies the formula (\ref{zwykla0}). To apply
Proposition \ref{stozek} we have to express the equivariant Chern
class $c^T(\1_{\P(X)})=1\otimes c_{SM}(\1_{\P(X)})$ by
$h=c_1^T(\O_{\P(X)}(1))$. The point is that the torus $T$ acts on
the fibers of the tautological bundle $\O(-1)$ with weight equal
to one, therefore the equivariant Chern class of
$\O_{\P^{n-1}}(1)$ is equal to
$$h=1\otimes x-t\otimes 1=x-t$$ under the identification
$$H_T^*(\P^{n-1})=\Q[\,t\,]\otimes H^*(\P^{n-1})=H^*(\P^{n-1})[\,t\,]\,.$$
Hence the equivariant Chern class of $\1_{\P(X)}$ can be written
as
$$c^T(\1_{\P(X)})=\Big(\sum_{i=0}^{n-1} a_i(h+t)^i\Big)\cap[\P^{n-1}]=
\Big(\sum_{i=0}^{n-1}\sum_{j=0}^i {i\choose j} a_i
t^{i-j}h^j\Big)\cap[\P^{n-1}]\,.$$
 Here the coefficient $b_0(t)$ of the expression (\ref{zdzial}) is equal
 to
 $$b_0(t)=\sum_{i=0}^{n-1} a_it^i\,.$$
By Proposition \ref{stozek} we obtain the claim.\qed
\s

We see that formally the Chern-MacPherson-Schwartz class of $\1_X$
in $\P^n$ and the equivariant Chern class in $\C^n$ satisfy the
same formula. The equivariant approach has the advantage that we
have for free the Chern class of the product
$$c^{T\times T}(\1_{X\times Y})=c^T(\1_X)\times c^T(\1_Y)\,.$$
Further we can restrict the Chern class of the product via
diagonal inclusion $T\hookrightarrow T\times T$ to obtain
 $c^T(\1_{X\times Y})$.
With the original approach the proof of the above property was a
bit demanding, see \cite[Lemma 3.10]{AMfe}

\section{Computing equivariant Chern classes without resolution of
singularities}\label{metoda}

Below we sketch a method of computing the equivariant Chern class
of a $T$-invariant singular variety not using a resolution of
singularities. The calculi will be done in equivariant cohomology
and we will omit the Poincar\'e  duality in the notation. \s

Assume that the fixed point set of the action of the torus on a
complex manifold $M$ is discrete. For a given class $x\in
H^k_T(M)$ of degree $k<2\,\dim(M)$ the integral $\int_{M}x$
vanishes. By the Localization Theorem also the sum $\sum_{p\in
M^T} \frac {x_{|p}}{e_p}$ has to vanish. In particular if
$x=c^T(\1_X)$, then except from the zero degree
\begin{equation}\sum_{p\in M^T}
\frac {c^T(\1_X)_{|p}}{e_p}=0\,.\label{BVforC}\end{equation} This
relation between local equivariant Chern classes allows  in many
cases to compute them inductively. Suppose
$M^T=\{p_0,p_1,\dots,p_N\} $ and assume that we know all local
equivariant Chern classes for $p_1,p_2,\dots,p_N$. Then
\begin{equation}c^T(\1_X)_{p_0}=-\sum_{i=1}^N \frac{e_{p_0}}{e_{p_i}}c^T(\1_X)_{p_i}\label{suma}\end{equation}
except from the zero degree. For Grassmannians the quotient
$\frac{e_{p_0}}{e_{p_i}}$ simplifies remarkably.\s

The zero component of the local equivariant Chern class is easy.
If $p\in X^T$ then by Theorem \ref{eu} this class is equal to the
Euler class at the point $p$
\begin{equation}(c^T(\1_X)_{(0)})_{p }=e_p \in H^{2\,\dim(M)}_T(pt).
\label{topo}\end{equation} In fact this statement is the crucial
point for computation. Any other equivariant characteristic class
satisfies the relation \ref{suma}. The condition fixing the zero
equivariant Chern class and vanishing for the degrees higher than
the dimension of the ambient space makes the equivariant Chern
class unique.

\pgfdeclareimage[height=7.48cm,width=12.25cm]{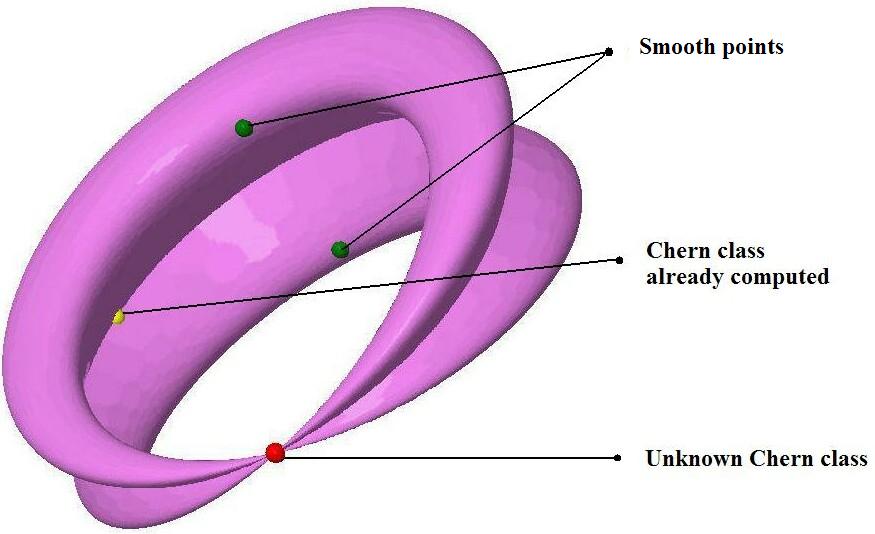}{strat2}
\begin{center} \pgfuseimage{strat2}\end{center}

\begin{center}{ Computation of the local equivariant Chern classes}\end{center}

Of course the inductive step of computation can be applied when
for a given singularity one can find a compact variety for which
this singularity is the only deepest one. If $X\subset \C^m$ is a
cone then taking the closure of $X$ in $\P^n$ will not introduce
new singularities. In the next section we present another
situation, when the compactifying variety is the Grassmannian.

\section{Computation of local equivariant Chern class of the determinant variety}\label{detvar}
\label{obliczenia}

Let us compute the local equivariant Chern class of the variety
$$\Omega_1^o(n)=\{\phi\in\Hom(\C^n,\C^n) \,:\,\det(\phi)=0\}\,.$$
Its  compactification in $Grass_n(\C^{2n})$ is the Schubert
variety of codimension one
$$\Omega_1(n)=\{W\,:\,W\cap\langle
\varepsilon_1,\varepsilon_2,\dots,\varepsilon_n\rangle\not=0\}\,.$$
We will apply the method sketched above. Let us start with $n=2$.
The canonical neighborhood of the point $p_{1,2}$ in
$\Grass_2(\C^4)$ is identified with
$$\Hom(span(\varepsilon_1,\varepsilon_2),span(\varepsilon_3,\varepsilon_4))$$
and the variety $\Omega_1(2)$ intersected with this neighbourhood
is exactly $\Omega^o_1(2)$. The corresponding elements of
$\Omega_1(2)$ are the planes spanned by the row-vectors of the
matrix
$$\left( \begin{matrix}1&0&a&b\cr 0& 1&c &d\end{matrix}\right)\,.$$
The equation of $\Omega_1(2)$ is
$$ \det\left( \begin{matrix}a&b\cr c
&d\end{matrix}\right)=0\,.$$ Before performing computations let us
draw the Goresky-Kottwitz-MacPherson graph (\cite[Th.~7.2]{GKM})
for $M=\Grass_2(\C^4)$ with the variety $\Omega_1(2)$ displayed.
\pgfdeclareimage[height=7.99cm,width=9cm]{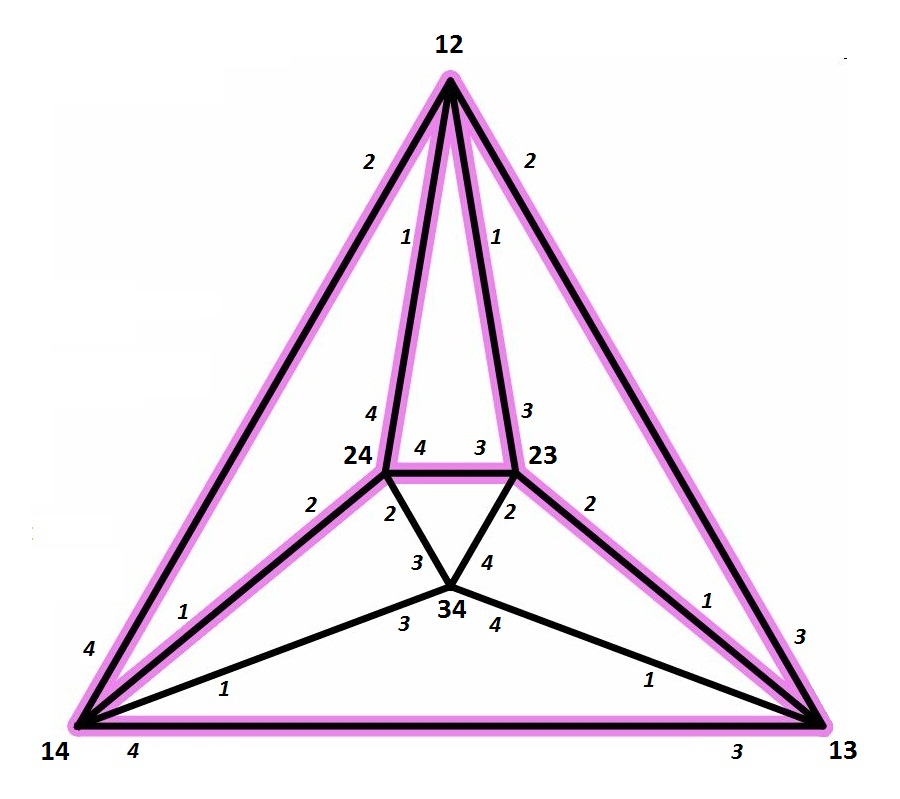}{v1n}
\begin{center} \pgfuseimage{v1n}\end{center}

\begin{center}{ Schubert variety $\Omega_1$ in $\Grass_2(\C^4)$.}\end{center}

The numbers attached to the edges indicate the weights of the $T$
actions along the one dimensional orbits. For example at the point
$p_{1,3}$ in the direction towards $p_{1,2}$ the action is by the
character $t_2-t_3$. The variety $\Omega_1(2)$ is singular at the
point $p_{1,2}$ and it is smooth at the remaining points. For
example at the point $p_{1,3}$ the coordinates are
$$\left( \begin{matrix}1&a&0&b\cr 0&c& 1&d\end{matrix}\right)$$
and the equation of $\Omega_1(2)$ is $b=0$. For that point the
local equivariant Chern class is equal to
$${ (t_4-t_1)}(1+t_2-t_1)(1+t_2-t_3)(1+t_4-t_3)\,.$$
The summand in the formula (\ref{BVforC}) is the following
$$\frac{{ (t_4-t_1)}(1+t_2-t_1)(1+t_2-t_3)(1+t_4-t_3)}{{ (t_4-t_1)}(t_2-t_1)(t_2-t_3)(t_4-t_3)}=$$
$$=\left(1+\frac{1}{t_2-t_1}\right)\left(1+\frac{1}{t_2-t_3}\right)\left(1+\frac{1}{t_4-t_3}\right)$$
We sum up the contribution coming from the fixed points $p_{1,3}$,
$p_{1,4}$, $p_{2,3}$, $p_{2,4}$, simplify and multiply by $-(t_3 -
t_1) (t_4 - t_1) (t_3 - t_2) (t_4 - t_2)$. We obtain \s $
(t_3+t_4-t_1-t_2)\hfill deg=1$ \s $(t_3+t_4-t_1-t_2)^2\hfill
deg=2$ \s $(t_3+t_4-t_1-t_2) (2 t_1 t_2-t_1 t_3-t_2 t_3-t_1
t_4-t_2 t_4+2 t_3 t_4)\hfill deg=3$ \s $-4 (t_3 - t_1) (t_4 - t_1)
(t_3 - t_2) (t_4 - t_2)\hfill deg=4$ \s \noindent The terms of
degree $<4$ coincide with the equivariant Chern class of
$\Omega_1(2)$ localized at the point $p_{1,2}$. The result is
symmetric in two groups of variables: $\{t_1,t_2\}$ and
$\{t_3,t_4\}$. The coefficients of the expansion in the basis of
the Schur functions
$$c^T(\1_{\Omega_1})=\sum a_{\rm I,J}S_{\rm I}(-t_1,-t_2)\cdot S_{\rm J}(t_3,t_4)$$
has the following coefficients:
$$
\begin{array}{c|cccccc}
& \it 0 &\it 1&\it 11& \it 2& \it 21& \it22\cr \hline \it 0 & & 1&
1& 1& 2& 1\cr \it 1 & 1 &1 & 3 &1& 1\cr \it 11 &1& 3 && 1\cr \it 2
& 1 &1& 1\cr \it 21 & 2& 1\cr \it 22 & 1\cr
\end{array}$$

Computations of the equivariant Chern class $\Omega_1(3)\subset
\Grass_3(\C^6)$ can be continued without problems by the same
method. At the points of the type $p_I$ with $|I\cap\{1,2,3\}|=1$
the variety is smooth, while at the points $p_I$ with
$|I\cap\{1,2,3\}|=2$ the singularity is of the type
$\Omega_1(2)_{p_{1,2}}\times \C^5$. We write the sum of fractions
according to the rule (\ref{suma}) and simplify. For example the
expression which has to be simplified to compute the degree one
is the following:
\begin{align}
&-\frac{(s_3-t_1) (s_3-t_2) (s_1-t_3) (s_2-t_3)}
{(s_3-s_1)(s_3-s_2) (t_1-t_3) (t_2-t_3)}\,(s_3-t_3)\;+\;sym.\;+\cr
&\frac{(s_3-t_1) (s_3-t_2) (s_1-t_3) (s_2-t_3)}
{(s_1-s_3)(s_2-s_3) (t_3-t_1)
(t_3-t_2)}\,(s_1+s_2-t_1-t_2)\;+\;sym. \label{k3}\end{align}
 (Here $s_1=t_4$, \;$s_2=t_5$, \;$s_3=t_6$.)
The given summands are the contributions coming from the points
$p_{3,4,5}$ and $p_{1,2,6}$. Of course the sum is equal to the
fundamental class
$$[\Omega_1]=s_1+s_2+s_3-t_1-t_2-t_3$$
(which may be computed in another way). This example shows how a
complicated rational functions may in fact lead to a simple
result. The difficulty lies in simplifying that expression. Higher
degree terms are much more complex. We write the final result in
the Schur basis
$$c^T(\1_{\Omega_1})=\sum a_{\rm I,J}S_{\rm I} (-t_1,-t_2,-t_3)\cdot S_{\rm J}(s_1,s_2,s_3)\,.$$
The coefficients are the following: {\footnotesize
$$
\begin{array}{c|cccccccccccccccccccc}
\!\!&\it\!\! 0 \!\!&\it\!\! 1 \!\!&\it\!\! 11 \!\!&\it\!\! 2
\!\!&\it\!\! 111 \!\!&\it\!\! 21 \!\!&\it\!\! 3 \!\!&\it\!\! 211
\!\!&\it\!\! 31 \!\!&\it\!\! 22 \!\!&\it\!\! 311 \!\!&\it\!\! 221
\!\!&\it\!\! 32 \!\!&\it\!\! \
321 \!\!&\it\!\! 222 \!\!&\it\!\! 33 \!\!&\it\!\! 331 \!\!&\it\!\! 322 \!\!&\it\!\! 332 \!\!&\it\!\! 333 \\
\hline \it 0 \!\!&\!\! \!\!&\!\! 1 \!\!&\!\! 2 \!\!&\!\! 2
\!\!&\!\! 4 \!\!&\!\! 5 \!\!&\!\! 1 \!\!&\!\! 9 \!\!&\!\! 3
\!\!&\!\! 4 \!\!&\!\! 6 \!\!&\!\! 9 \!\!&\!\! 3 \!\!&\!\! 8
\!\!&\!\! 4 \!\!&\!\! 1 \!\!&\!\! \ 3 \!\!&\!\! 6 \!\!&\!\! 3
\!\!&\!\! 1 \\ \it 1 \!\!&\!\! 1 \!\!&\!\! 4 \!\!&\!\! 8 \!\!&\!\!
5 \!\!&\!\! 12 \!\!&\!\! 12 \!\!&\!\! 2 \!\!&\!\! 19 \!\!&\!\! 5
\!\!&\!\! 8 \!\!&\!\! 8 \!\!&\!\! 16 \!\!&\!\! 4 \!\!&\!\! 8
\!\!&\!\! 10 \!\!&\!\! \ 1 \!\!&\!\! 2 \!\!&\!\! 4 \!\!&\!\! 1
\!\!&\!\! \\ \it 11 \!\!&\!\! 2 \!\!&\!\! 8 \!\!&\!\! 12 \!\!&\!\!
9 \!\!&\!\! 16 \!\!&\!\! 16 \!\!&\!\! 3 \!\!&\!\! 18 \!\!&\!\! 6
\!\!&\!\! 8 \!\!&\!\! 6 \!\!&\!\! 10 \!\!&\!\! 4 \!\!&\!\! 4
\!\!&\!\! \ \!\!&\!\! 1 \!\!&\!\! 1 \!\!&\!\! \!\!&\!\! \!\!&\!\!
\\ \it 2 \!\!&\!\! 2 \!\!&\!\! 5 \!\!&\!\! 9 \!\!&\!\! 4 \!\!&\!\!
11 \!\!&\!\! 9 \!\!&\!\! 1 \!\!&\!\! 13 \!\!&\!\! 2 \!\!&\!\! 5
\!\!&\!\! 3 \!\!&\!\! 10 \!\!&\!\! 1 \!\!&\!\! 2 \!\!&\!\! 5
\!\!&\!\! \ \!\!&\!\! \!\!&\!\! 1 \!\!&\!\! \!\!&\!\! \\ \it 111
\!\!&\!\! 4 \!\!&\!\! 12 \!\!&\!\! 16 \!\!&\!\! 11 \!\!&\!\! 8
\!\!&\!\! 16 \!\!&\!\! 4 \!\!&\!\! \!\!&\!\! 6 \!\!&\!\! 10
\!\!&\!\! \!\!&\!\! \!\!&\!\! 4 \!\!&\!\! \!\!&\!\! \ \!\!&\!\! 1
\!\!&\!\! \!\!&\!\! \!\!&\!\! \!\!&\!\! \\ \it 21 \!\!&\!\! 5
\!\!&\!\! 12 \!\!&\!\! 16 \!\!&\!\! 9 \!\!&\!\! 16 \!\!&\!\! 14
\!\!&\!\! 2 \!\!&\!\! 15 \!\!&\!\! 3 \!\!&\!\! 5 \!\!&\!\! 3
\!\!&\!\! 5 \!\!&\!\! 1 \!\!&\!\! 1 \!\!&\!\! \ \!\!&\!\!
\!\!&\!\! \!\!&\!\! \!\!&\!\! \!\!&\!\! \\ \it 3 \!\!&\!\! 1
\!\!&\!\! 2 \!\!&\!\! 3 \!\!&\!\! 1 \!\!&\!\! 4 \!\!&\!\! 2
\!\!&\!\! \!\!&\!\! 3 \!\!&\!\! \!\!&\!\! 1 \!\!&\!\! \!\!&\!\! 2
\!\!&\!\! \!\!&\!\! \!\!&\!\! 1 \!\!&\!\! \!\!&\!\! \ \!\!&\!\!
\!\!&\!\! \!\!&\!\! \\ \it 211 \!\!&\!\! 9 \!\!&\!\! 19 \!\!&\!\!
18 \!\!&\!\! 13 \!\!&\!\! \!\!&\!\! 15 \!\!&\!\! 3 \!\!&\!\!
\!\!&\!\! 3 \!\!&\!\! 5 \!\!&\!\! \!\!&\!\! \!\!&\!\! 1 \!\!&\!\!
\!\!&\!\! \ \!\!&\!\! \!\!&\!\! \!\!&\!\! \!\!&\!\! \!\!&\!\! \\
\it 31 \!\!&\!\! 3 \!\!&\!\! 5 \!\!&\!\! 6 \!\!&\!\! 2 \!\!&\!\! 6
\!\!&\!\! 3 \!\!&\!\! \!\!&\!\! 3 \!\!&\!\! \!\!&\!\! 1 \!\!&\!\!
\!\!&\!\! 1 \!\!&\!\! \!\!&\!\! \!\!&\!\! \!\!&\!\! \!\!&\!\! \
\!\!&\!\! \!\!&\!\! \!\!&\!\! \\ \it 22 \!\!&\!\! 4 \!\!&\!\! 8
\!\!&\!\! 8 \!\!&\!\! 5 \!\!&\!\! 10 \!\!&\!\! 5 \!\!&\!\! 1
\!\!&\!\! 5 \!\!&\!\! 1 \!\!&\!\! \!\!&\!\! 1 \!\!&\!\! \!\!&\!\!
\!\!&\!\! \!\!&\!\! \!\!&\!\! \ \!\!&\!\! \!\!&\!\! \!\!&\!\!
\!\!&\!\! \\ \it 311 \!\!&\!\! 6 \!\!&\!\! 8 \!\!&\!\! 6 \!\!&\!\!
3 \!\!&\!\! \!\!&\!\! 3 \!\!&\!\! \!\!&\!\! \!\!&\!\! \!\!&\!\! 1
\!\!&\!\! \!\!&\!\! \!\!&\!\! \!\!&\!\! \!\!&\!\! \!\!&\!\! \
\!\!&\!\! \!\!&\!\! \!\!&\!\! \!\!&\!\! \\ \it 221 \!\!&\!\! 9
\!\!&\!\! 16 \!\!&\!\! 10 \!\!&\!\! 10 \!\!&\!\! \!\!&\!\! 5
\!\!&\!\! 2 \!\!&\!\! \!\!&\!\! 1 \!\!&\!\! \!\!&\!\! \!\!&\!\!
\!\!&\!\! \!\!&\!\! \!\!&\!\! \!\!&\!\! \ \!\!&\!\! \!\!&\!\!
\!\!&\!\! \!\!&\!\! \\ \it 32 \!\!&\!\! 3 \!\!&\!\! 4 \!\!&\!\! 4
\!\!&\!\! 1 \!\!&\!\! 4 \!\!&\!\! 1 \!\!&\!\! \!\!&\!\! 1
\!\!&\!\! \!\!&\!\! \!\!&\!\! \!\!&\!\! \!\!&\!\! \!\!&\!\!
\!\!&\!\! \!\!&\!\! \!\!&\!\! \ \!\!&\!\! \!\!&\!\! \!\!&\!\! \\
\it 321 \!\!&\!\! 8 \!\!&\!\! 8 \!\!&\!\! 4 \!\!&\!\! 2 \!\!&\!\!
\!\!&\!\! 1 \!\!&\!\! \!\!&\!\! \!\!&\!\! \!\!&\!\! \!\!&\!\!
\!\!&\!\! \!\!&\!\! \!\!&\!\! \!\!&\!\! \!\!&\!\! \ \!\!&\!\!
\!\!&\!\! \!\!&\!\! \!\!&\!\! \\ \it 222 \!\!&\!\! 4 \!\!&\!\! 10
\!\!&\!\! \!\!&\!\! 5 \!\!&\!\! \!\!&\!\! \!\!&\!\! 1 \!\!&\!\!
\!\!&\!\! \!\!&\!\! \!\!&\!\! \!\!&\!\! \!\!&\!\! \!\!&\!\!
\!\!&\!\! \!\!&\!\! \ \!\!&\!\! \!\!&\!\! \!\!&\!\! \!\!&\!\! \\
\it 33 \!\!&\!\! 1 \!\!&\!\! 1 \!\!&\!\! 1 \!\!&\!\! \!\!&\!\! 1
\!\!&\!\! \!\!&\!\! \!\!&\!\! \!\!&\!\! \!\!&\!\! \!\!&\!\!
\!\!&\!\! \!\!&\!\! \!\!&\!\! \!\!&\!\! \!\!&\!\! \!\!&\!\! \
\!\!&\!\! \!\!&\!\! \!\!&\!\! \\ \it 331 \!\!&\!\! 3 \!\!&\!\! 2
\!\!&\!\! 1 \!\!&\!\! \!\!&\!\! \!\!&\!\! \!\!&\!\! \!\!&\!\!
\!\!&\!\! \!\!&\!\! \!\!&\!\! \!\!&\!\! \!\!&\!\! \!\!&\!\!
\!\!&\!\! \!\!&\!\! \ \!\!&\!\! \!\!&\!\! \!\!&\!\! \!\!&\!\! \\
\it 322 \!\!&\!\! 6 \!\!&\!\! 4 \!\!&\!\! \!\!&\!\! 1 \!\!&\!\!
\!\!&\!\! \!\!&\!\! \!\!&\!\! \!\!&\!\! \!\!&\!\! \!\!&\!\!
\!\!&\!\! \!\!&\!\! \!\!&\!\! \!\!&\!\! \!\!&\!\! \ \!\!&\!\!
\!\!&\!\! \!\!&\!\! \!\!&\!\! \\ \it 332 \!\!&\!\! 3 \!\!&\!\! 1
\!\!&\!\! \!\!&\!\! \!\!&\!\! \!\!&\!\! \!\!&\!\! \!\!&\!\!
\!\!&\!\! \!\!&\!\! \!\!&\!\! \!\!&\!\! \!\!&\!\! \!\!&\!\!
\!\!&\!\! \!\!&\!\! \ \!\!&\!\! \!\!&\!\! \!\!&\!\! \!\!&\!\! \\
\it 333 \!\!&\!\! 1 \!\!&\!\! \!\!&\!\! \!\!&\!\! \!\!&\!\!
\!\!&\!\! \!\!&\!\! \!\!&\!\! \!\!&\!\! \!\!&\!\! \!\!&\!\!
\!\!&\!\! \!\!&\!\! \!\!&\!\! \!\!&\!\! \!\!&\!\! \ \!\!&\!\!
\!\!&\!\! \!\!&\!\! \!\!&\!\!
\end{array}
$$
}

\noindent We note that all the coefficients are nonnegative.

While computing the equivariant Chern class of $\Omega_1(4)\subset
\Grass_4(\C^8)$ appears a problem with the size of the
expressions, since $\dim(\Grass_4(\C^8))=16$ and $\dim(T)=8$. In a
polynomial of degree 15 in 8 variables there are $$490\,314\;{\rm
monomials.}$$ The expression is a sums of 68 fractions with
factors $t_i-t_j$ in denominators. We might have used another
compactification of $\C^{16}$, e.g.~the projective space $
\P^{16}$. There are less fixed points, but the denominators are
more complicated. They are of the form $\prod
[(t_i-t_j)-(t_k-t_\ell)]$.

One practical solution appears naturally. The fixed points can be
divided into groups with $|I\cap\{1,2,\dots,n\}|$ fixed. Let
$f_k(u_\bullet,v_\bullet)$ be the expression for the local
equivariant Chern class of $\Omega_1(k)$ with
$u_\bullet=(t_1,t_2,\dots,t_k)$ and
$v_\bullet=(t_{k+1},t_{k+2},\dots,t_{2k})$. The local equivariant
Chern $c^T(\1_{\Omega_1(n)})_{|p_{1,2,\dots,n}}$ class can by
computed by the formula (\ref{suma}), which becomes
\begin{equation}
-\sum_{k=1}^{n-1}\sum_{\begin{matrix}I\subset\{1,2,\dots,2n\}\cr|I|=n,\,|I\cap\{1,2,\dots,n\}|=k\end{matrix}}
\frac{e_{p_{1,2,\dots,n}}}{e_{p_I}}\,f_k(I)\,g_k(I)\,,\label{grgr}\end{equation}
where $f_k(I)$ depends on the two group of variables
$$u_\bullet=t_{I\cap\{1,2,\dots,n\}}\quad{\rm and}\quad v_\bullet=t_{\{n+1,n+2,\dots,2n\}\setminus I}$$
and $g_k(I)$ is the equivariant Chern class of the singular
stratum of the type $\Omega_1(k)$. The factors in the quotients
$\frac{e_{p_{1,2,\dots,n}}}{e_{p_I}}$ cancel out partially and
miraculously all the summands for a fixed $k$ turn out to be
integral. For $n=3$ and degree one the summands are given by the
formula (\ref{k3}). \s Such a division of fixed points has a
geometric meaning. In fact we deal with the partial localization
(see \S\ref{par-loc}). Consider the action of the subtorus $\C^*$
acting on $\C^{2n}$ with weight 1 on the first $n$ coordinates and
with the weight $-1$ on the remaining coordinates. Then the fixed
point set decomposes into disjoint union of the products of the
Grassmannians:
$$\Grass_n(\C^{2n})^{\C^*}= \bigsqcup_{k=0}^n \Grass_{k}(\C^n)\times \Grass_{n-k}(\C^n).$$
The summand for $k=0$ consists of one point $$\{0\}\oplus\langle
\varepsilon_{n+1},\varepsilon_{n+2},\dots,\varepsilon_{2n}\rangle\,,$$
which does not belong to $\Omega_1(n)$, while for $k=n$ we have
$$\langle \varepsilon_1,\varepsilon_2,\dots,\varepsilon_n\rangle\oplus\{0\}\,,$$ the point which we are concerned with.
Let ${\cal R}_{k}$ and ${\cal Q}_{k}$ be the tautological and the
quotient bundles over $\Grass_k(\C^n)$. The result of the sum
(\ref{grgr}) is equal to
$$-(-1)^{n-k}\,\sum_{k=1}^{n-1}\int_{\Grass_{k}(\C^n)\times \Grass_{n-k}(\C^n)}
\left[f_{k}({\cal R}_{k},{\cal Q}_{n-k})\cdot \bar g_k\right]\,,$$
where $\bar g_k$ is (up to multiplication by a certain Euler
class) the equivariant Chern class of the stratum of the
singularity type $\Omega_1(k)$. Precisely\s $\bar g_k= e({\cal
Q}^*_{k}\otimes{\cal Q}_{n-k})\cdot e({\cal R}^*_{k}\otimes{\cal
R}_{n-k})\cdot$\s \hfill$ \cdot c({\cal R}^*_{k}\otimes{\cal
Q}_{k})\cdot c({\cal R}^*_{n-k}\otimes{\cal Q}_{n-k})\cdot c({\cal
Q}_{k}\otimes{\cal R}^*_{n-k})\,.$\s \noindent Using Fubini
theorem we do not have to simplify a large expression in one step
and we arrive to the result relatively quickly. Also knowing the
Schur expansion of the functions $f_k$ one can apply Theorem
\ref{jlp}.

\section{GKM-relations}

Less time-consuming method of computation of the local equivariant
Chern class is based on the relation discovered by Chen-Skjelbred
\cite{CS}, called GKM-relations after the rediscovery in
\cite{GKM}. These relation allow us to determine the local
equivariant Chern class at the point $p_I$ knowing only the local
equivariant Chern classes at the neighbouring points in the
GKM-graph. This is so since
$$c^T(\1_V)_{|p_{\rm I}}\equiv c^T(\1_V)_{|p_{\rm J}}\; {\mathop{\rm modulo}}\; (t_i-t_j)\,,$$
whenever $${\rm J}=({\rm I}-\{i\})\cup\{j\}\,.$$ Again the method
works for all the degrees smaller then the dimension of the
Grassmannian, since the intersection of the ideals $(t_i-t_j)$ is
contained in the degree greater or equal to the dimension of
the Grassmannian. That is so for any GKM-space. Now the problem of
simplifying huge rational function is replaced by solving a
relatively small system of linear equations.

\section{The result for $\Grass_4(\C^8)$}\label{monstrum}

Let us write the local equivariant Chern class in the Schur basis
$$c^T(\1_{\Omega_1(4)})_{|p_{1,2,3,4}}=\sum a_{I,J}S_{\rm
I}(-t_1,-t_2,-t_3,-t_4)\cdot S_{\rm J}(t_5,t_6,t_7,t_8)\,.$$ Just
to quench readers curiosity we show here the most interesting
fragment of the table of coefficients.\s

It is hard not to have impression that there should be a way of
writing down this equivariant Chern class in a compact way. For
example the equivariant Chern class of the tangent bundle written
in the Schur basis is as much complicated as ours, but it is just
$c(\Hom({\cal R}_n,{\cal Q}_n))$.\s

\pgfdeclareimage[height=21cm,width=14cm]{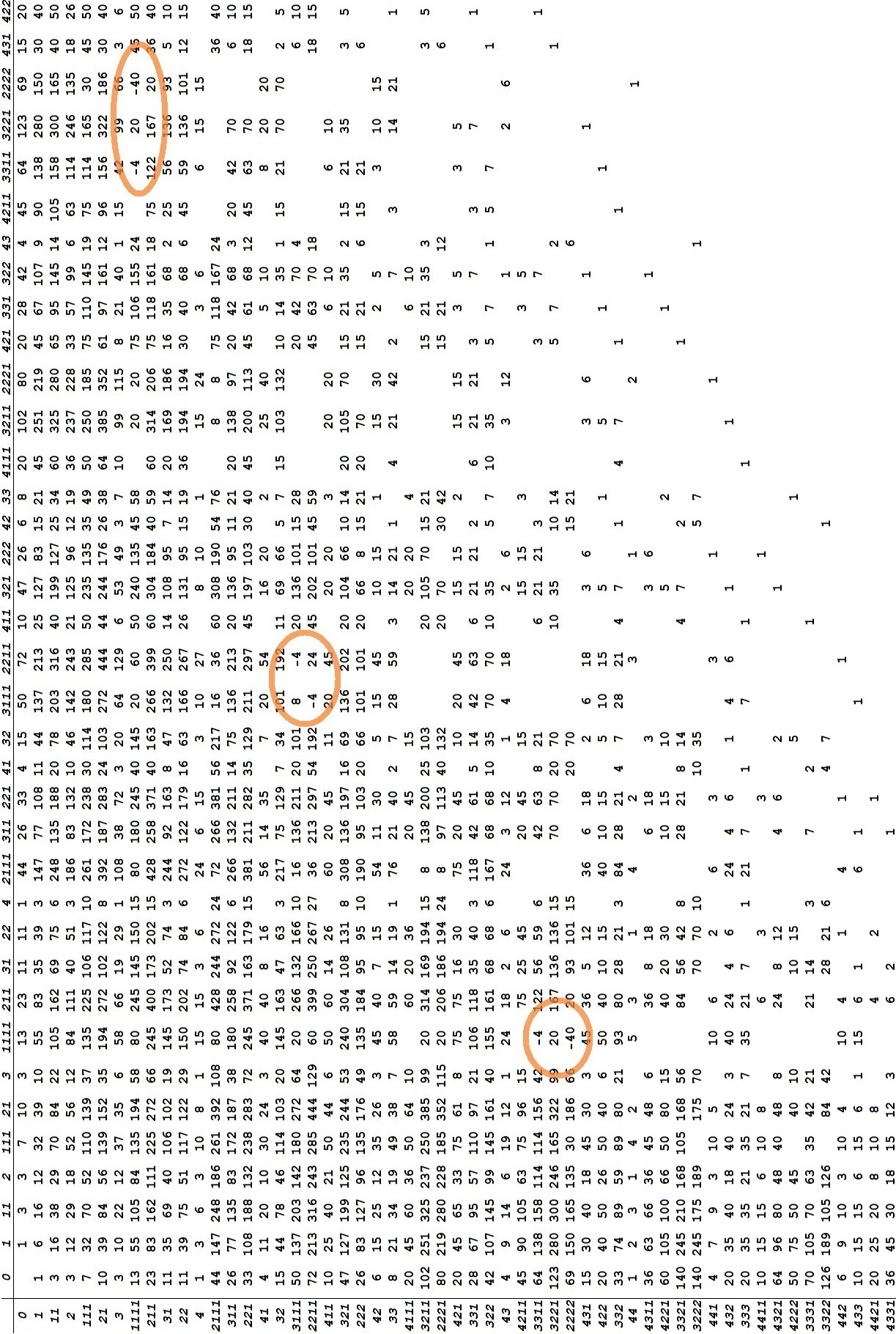}{k4k}

\noindent\hskip -20pt \pgfuseimage{k4k}

It turns out that the local equivariant Chern class of
$\Omega_1(4)$ is a positive combination of monomials in
$-t_1,-t_2,-t_3,-t_4,t_5,t_6,t_7,t_8$. As one can see it is {\it
not} a positive combinations of products of Schur functions.
Fortunately we do not have a contradiction with the conjecture of
Aluffi and Mihalcea \cite{AlMi} which says that the
Chern-Schwartz-MacPherson  classes are effective. Note that the
Schubert varieties are only $T$-invariant, and the Theorem
\ref{pos0} does not apply. Instead we have a freedom with choosing
the basis of weights. The local equivariant Chern class is a
polynomial in $u_{i,j}=t_i-t_j$. To write $c^T(\1_X)$ in a unique
way we chose a spanning tree of the full graph with vertices
$1,2,\dots,2n$. The edge between $i$ and $j$ (with the orientation
forced by the partition) corresponds to the generator $t_j-t_i$.
Some choices lead to an expression with nonnegative coefficients.

\s

\hfil {Positive monomial bases for $\Grass_2(\C^4)$}

\pgfdeclareimage[height=1.72cm,width=11.82cm]{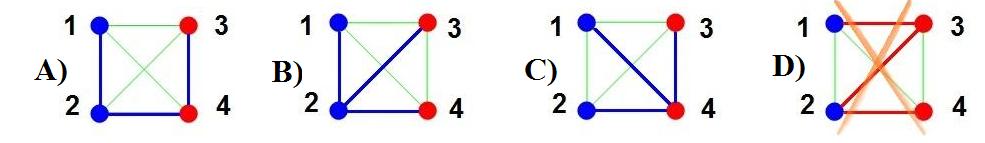}{kolumnap}
\begin{center} \pgfuseimage{kolumnap}\end{center}
\vskip-10pt

A) $t_2-t_1,\quad t_4-t_2,\quad t_3-t_4$

B) $t_2-t_1,\quad t_3-t_2,\quad t_4-t_2$

C) $t_4-t_1,\quad t_4-t_2,\quad t_4-t_3$

D) $t_3-t_1,\quad t_3-t_2,\quad t_4-t_2$~~~~ this is not a
positive basis \s

\noindent The positivity condition for a graph is the
following:\begin{itemize}\item Characters of the tangent
representation are nonnegative sums of base elements.\end{itemize}
That in fact supports the conjecture of Aluffi and Mihalcea in a
stronger, equivariant version. \s The original, nonequivariant
version was checked by B.~Jones \cite{Jon} for cells in
$\Grass_m(\C^n)$ for $m\leq 3$. In his computations equivariant
cohomology and the Localization Theorem was used to compute the
push-forward of classes from a resolution.

\section{Further directions of work}
Several goals have not been reached so far. The most obvious
directions of further work would be:
\begin{itemize}
\item deduce positivity results,
\item study global equivariant Chern classes of Schubert varieties
and open cells,
\item in particular relate our computations to the determinant
formulas of \cite{AlMi} and the combinatorial interpretation of
\cite{Mi},
\item develop a suitable calculus of symmetric rational functions to handle expressions appearing in
the Berline-Vergne formula for Grassmannians.
\end{itemize}
We hope to realize this program in future.

\end{document}